\documentclass[11pt]{amsart}

\usepackage{amsmath}
\usepackage{amssymb}

\newtheorem{theorem}{Theorem}[section]
\newtheorem{claim}[theorem]{Claim}

\theoremstyle{definition}
\newtheorem{definition}[theorem]{Definition}

\newtheorem{question}[theorem]{Question}

\theoremstyle{remark}
\newtheorem{remark}[theorem]{Remark}

\newtheorem{notation}[theorem]{Notation}

\newcount\skewfactor
\def\mathunderaccent#1#2 {\let\theaccent#1\skewfactor#2
\mathpalette\putaccentunder}
\def\putaccentunder#1#2{\oalign{$#1#2$\crcr\hidewidth
\vbox to.2ex{\hbox{$#1\skew\skewfactor\theaccent{}$}\vss}\hidewidth}}
\def\name{\mathunderaccent\tilde-3 }

\def\smallbox#1{\leavevmode\thinspace\hbox{\vrule\vtop{\vbox
   {\hrule\kern1pt\hbox{\vphantom{\tt/}\thinspace{\tt#1}\thinspace}}
   \kern1pt\hrule}\vrule}\thinspace}

\newcommand{\cf}{{\rm cf}}

\def\qedref#1{$\qed_{\reforiginal{#1}}$}

\title{Amenable colorings}
\author{Shimon Garti}
\address{Institute of Mathematics,
 The Hebrew University of Jerusalem,
 Jerusalem 91904, Israel}
\email{shimon.garty@mail.huji.ac.il}

\begin{document}
\let\labeloriginal\label
\let\reforiginal\ref

\begin{abstract}
Let $\kappa$ be any regular cardinal.
Assuming the existence of a huge cardinal above $\kappa$, we prove the consistency of $\binom{\kappa^{++}}{\kappa^+}\rightarrow \binom{\tau}{\kappa^+}^{1,1}_\kappa$ for every ordinal $\tau<\kappa^{++}$.
Likewise, we prove that $\binom{\aleph_2}{\aleph_1}\rightarrow_{\mathcal{A}} \binom{\aleph_2}{\aleph_1}^{1,1}_2$ is consistent when $\mathcal{A}$ is strongly closed under countable intersections.
\end{abstract}

\maketitle

\newpage

\section{Introduction}

The strong polarized partition relation $\binom{\lambda}{\kappa}\rightarrow \binom{\lambda}{\kappa}^{1,1}_\theta$ means that for every coloring $c:\lambda\times\kappa\rightarrow\theta$ there exists a monochromatic product $A\times B$ so that $|A|=\lambda$ and $|B|=\kappa$. 
Two major problems stand in the front. For any given infinite cardinal $\kappa$ we ask about the pair $(\lambda, \kappa)$ with respect to $\lambda=\kappa^+$ and $\lambda=2^\kappa$. Outside the interval $[\kappa^+,2^\kappa]$ the question becomes uninteresting. Let us try to explain why.

Firstly, we may always assume that $\kappa\leq\lambda$, since the notation of the strong polarized relation is symmetric. A simple coloring shows that $\binom{\kappa}{\kappa}\nrightarrow \binom{\kappa}{\kappa}^{1,1}_\theta$ for every $\kappa$ (actually, a stronger negation can be proved), so our investigation begins with $\lambda\geq\kappa^+$. The interval $[\kappa^+,2^\kappa]$ exhibits non-trivial demeanor, as positive and negative statements can be proved both for specific $\lambda\in[\kappa^+,2^\kappa]$ and for the behavior of the entire interval. If $\lambda=\cf(\lambda)>2^\kappa$ then $\binom{\lambda}{\kappa}\rightarrow \binom{\lambda}{\kappa}^{1,1}_\theta$ follows from the fact that $|\mathcal{P}(\kappa)|=2^\kappa$, so the right-hand component of every coloring will be the same for $\lambda$-many ordinals. If $\lambda>2^\kappa$ is a singular cardinal then the behavior of $\lambda$ with respect to the relation $\binom{\lambda}{\kappa}\rightarrow \binom{\lambda}{\kappa}^{1,1}_\theta$ is determined by the behavior of $\cf(\lambda)$ with respect to the relation $\binom{\cf(\lambda)}{\kappa}\rightarrow \binom{\cf(\lambda)}{\kappa}^{1,1}_\theta$. Hence a knowledge of the pertinent relations for $\lambda\in[\kappa^+,2^\kappa]$ gives a full knowledge for every $\lambda$. As a reference to the facts mentioned in this paragraph we suggest Chapter 4 in \cite{williams} (in particular, Theorem 4.14 and Lemma 4.2.7).

In this paper we focus on the pair $(\kappa^+,\kappa)$. A negative consistency relation can be forced for every $\kappa$ since $2^\kappa=\kappa^+$ implies $\binom{\kappa^+}{\kappa}\nrightarrow \binom{\kappa^+}{\kappa}^{1,1}_2$. For a positive consistency relation it seems natural to classify infinite cardinals into three categories.

If $\kappa$ is a large cardinal (including the case $\kappa=\aleph_0$) then one can force $\binom{\kappa^+}{\kappa}\rightarrow \binom{\kappa^+}{\kappa}^{1,1}_2$ by increasing the splitting number $\mathfrak{s}_\kappa$. 
Assuming that $\kappa=\cf(\kappa)$ it is known that $\mathfrak{s}_\kappa\geq\kappa$ iff $\kappa$ is strongly inaccessible, and $\mathfrak{s}_\kappa>\kappa$ iff $\kappa$ is weakly compact. So if one wishes to force $\binom{\kappa^+}{\kappa}\rightarrow \binom{\kappa^+}{\kappa}^{1,1}_2$ by increasing the splitting number, at least weak compactness must be assumed. Moreover, one needs $\mathfrak{s}_\kappa>\kappa^+$, and it is unknown if this setting is possible for mild large cardinals.
It has been done for every supercompact cardinal (see \cite{MR3000439} and \cite{MR3201820}), and recently also when $\kappa$ is measurable with large enough Mitchell order (see \cite{MR3436372}). 
Of course, perhaps $\binom{\kappa^+}{\kappa}\rightarrow \binom{\kappa^+}{\kappa}^{1,1}_2$ can be forced without increasing $\mathfrak{s}_\kappa$, so it is still open for small large cardinals whether this strong relation is forceable (see \cite{1012}, Question 4.4).

The second category is singular cardinals.
If $\kappa$ is a singular cardinal then we have a comprehensive answer, as $\binom{\kappa^+}{\kappa}\rightarrow \binom{\kappa^+}{\kappa}^{1,1}_2$ can be forced at every singular cardinal (see \cite{MR2987137} and \cite{1012}).

The third category is successor cardinals. It is unknown whether $\binom{\kappa^+}{\kappa}\rightarrow \binom{\kappa^+}{\kappa}^{1,1}_2$ can be forced on such cardinals.
In order to deal with this case we consider two different directions. The first one is based on the concept of amenable colorings, and the second is related to the concept of almost strong relations. Let us explain, shortly, the main idea of these concepts.

Given a collection $\mathcal{A}\subseteq[\kappa]^\kappa$ we focus on a coloring $c:\lambda\times\kappa\rightarrow\theta$ such that every fiber $\{\gamma\}\times\kappa$ has a monochromatic subset of the form $\{\gamma\}\times A_\gamma$ for some $A_\gamma\in\mathcal{A}$. Notice that the usual polarized relation is just the special case of $\mathcal{A}=[\kappa]^\kappa$.
In the next section we shall focus on the pair $(\aleph_2,\aleph_1)$ for which we shall prove the consistency of $\binom{\aleph_2}{\aleph_1}\rightarrow_{\mathcal{A}} \binom{\aleph_2}{\aleph_1}^{1,1}_2$ with respect to a suitable $\mathcal{A}$. The precise definitions and required properties are given at the beginning of this section, but the main point is that amenability may give consistency results even with full monochromatic products.

In the last section we concentrate on the common polarized relation, but our monochromatic product is just \emph{almost strong}. For colorings defined on $\lambda\times\kappa$ it means that the left-hand component can be of order type $\tau$ for every $\tau<\lambda$. Again, the precise defintion will be given at the beginning of the last section, but the theorem reads as follows:
The relation $\binom{\kappa^{++}}{\kappa^+}\rightarrow \binom{\tau}{\kappa^+}^{1,1}_\kappa$ for every ordinal $\tau<\kappa^{++}$ can be forced at every regular cardinal $\kappa$ (by assuming the presence of a huge cardinal above $\kappa$ in the ground model).

We use standard notation. If $A,B\subseteq\kappa$ then $A\subseteq^*B$ iff $|A\setminus B|<\kappa$. 
If $\kappa=\cf(\kappa)<\lambda$ then $S_\kappa^\lambda=\{\alpha<\lambda:\cf(\alpha)=\kappa\}$. Notice that $S^\lambda_\kappa$ is a stationary subset of $\lambda$. We use the Jerusalem forcing notation, i.e. $p\leq q$ means that the condition $q$ is stronger than $p$.
A forcing notion $\mathbb{P}$ is $\kappa$-centered iff $\mathbb{P}$ can be decomposed into $\kappa$-many subsets, each of which consists of pairwise compatible conditions. If $\mathcal{A}\subseteq[\kappa]^\kappa$ then $\mathcal{A}$ is strongly closed under intersections iff the cardinality of $a\cap b$ is $\kappa$ for every $a,b\in\mathcal{A}$. Similarly define the notion of $\mathcal{A}$ being strongly closed under countable intersections, and so on.

Several generalizations of Martin's Axiom for $\aleph_1$ are known in the literature. We shall make use of Shelah's version (but the variants of Baumgartner and Laver can serve as well):

\begin{theorem}
\label{sssshelah} Martin's Axiom for $\aleph_1$. \newline
One can force $2^{\aleph_0}=\aleph_1\wedge 2^{\aleph_1}>\aleph_2$, and if $\mathbb{P}$ is a forcing notion of size less than $2^{\aleph_1}$ satisfying the following three requirements:
\begin{enumerate}
\item [$(a)$] Each pair of compatible conditions has a least upper bound in $\mathbb{P}$.
\item [$(b)$] Every countable increasing sequence of conditions has a least upper bound in $\mathbb{P}$.
\item [$(c)$] If $\{p_i:i<\aleph_2\}\subseteq\mathbb{P}$ then there is a club $C\subseteq\aleph_2$ and a regressive function $f:\aleph_2\rightarrow\aleph_2$ so that for $\alpha,\beta\in C\cap S^{\aleph_2}_{\aleph_1}$ if $f(\alpha)=f(\beta)$ then $p_\alpha\parallel p_\beta$.
\end{enumerate}
then there is a generic filter $G\subseteq\mathbb{P}$ which intersects any given collection of $\kappa$-many dense subsets, when $\kappa<2^{\aleph_1}$.
\end{theorem}

We shall refer to the above statement as the generalized Martin's axiom.
The proof of the theorem appears in \cite{MR0505492}. We indicate that if $\kappa$ satisfies $\alpha<\kappa\Rightarrow\alpha^{\aleph_0}<\kappa$ then the assumption $|\mathbb{P}|<2^{\aleph_1}$ can be omitted (as shown in the above mentioned paper). Observe also that if $\mathbb{P}$ is $\aleph_1$-centered then requirement $(c)$ follows.

A cardinal $\kappa$ is huge iff there exists an elementary embedding $\jmath:{\rm V}\rightarrow M$ so that $\kappa={\rm crit}(\jmath)$ and ${}^{\jmath(\kappa)}M\subseteq M$. An ideal $\mathcal{I}$ is $(\mu,\mu,\theta)$-saturated iff for every collection $\mathcal{A}=\{A_\alpha: \alpha<\mu\}\subseteq\mathcal{I}^+$ there exists a sub-collection $\mathcal{B}\in[\mathcal{A}]^\mu$ such that $\mathcal{C}\in[\mathcal{B}]^\theta \Rightarrow \bigcap\limits_{\alpha\in\mathcal{C}}A_\alpha\in\mathcal{I}^+$. The following theorem belongs to Laver, \cite{MR673792}:

\begin{theorem}
\label{lavthm} Assume there exists a huge cardinal, and $\theta$ is a regular cardinal below this huge cardinal. \newline 
Then it is consistent that there is a $\theta^+$-complete and even normal ideal $\mathcal{I}$ over $\theta^+$ which is $(\theta^{++},\theta^{++},\theta)$-saturated. The existence of such an ideal can be forced also with $2^\theta=\theta^+$, and it preserves cardinalities and cofinalities in the interval $[\aleph_1,\theta]$.
\end{theorem}

\hfill \qedref{lavthm}

The idea behind the proof of the theorem is captured in the words of Prince Humperdinck: ``Someone has beaten a giant" (\cite{pbride}, p. 191). By collapsing a huge cardinal one can preserve some of its qualities, resulting in the existence of a sufficiently saturated ideal.
By and large, good combinatorial theorems hold over large cardinals, since the existence of a complete ultrafilter gives large monochromatic sets. However, a saturated ideal can play the r\^ole of an ultrafilter under suitable circumstances.

I wish to thank the referee of the paper for an extraordinary work, including both mathematical corrections and meaningful improvements of the presentation. This includes an elegant argument which simplified the proof of Theorem \ref{galthm}. I also thank Yair Hayut for his help.

\newpage 

\section{Amenability}

We begin with the concept of amenability:

\begin{definition}
\label{aaaa} Amenable coloring. \newline
Let $c:\lambda\times\kappa\rightarrow\theta$ be a coloring, and assume $\mathcal{A}\subseteq\mathcal{P}(\kappa)$. \newline 
We say that $c$ is $\mathcal{A}$-amenable if for every $\gamma<\lambda$ there are $i_\gamma<\theta$ and $A_\gamma\in\mathcal{A}$ so that $\delta\in A_\gamma \Rightarrow c(\gamma,\delta)=i_\gamma$.
\end{definition}

With the above definition we introduce the following notation:

\begin{notation}
\label{aanotation} $\rightarrow_{\mathcal{A}}$. \newline 
We say that $\binom{\lambda}{\kappa}\rightarrow_{\mathcal{A}} \binom{\lambda}{\kappa}^{1,1}_\theta$ holds iff for every $c:\lambda\times\kappa \rightarrow\theta$ which is $\mathcal{A}$-amenable there are $A\in[\lambda]^\lambda, B\in[\kappa]^\kappa$ and a color $\iota<\theta$ so that $c\upharpoonright(A\times B)=\{\iota\}$.
\end{notation}

The main theorem of this section establishes a positive consistency result of the strong relation for suitable amenability. In order to motivate the positive direction, we introduce the following:

\begin{claim}
\label{ggggch} Negative relations and GCH. \newline 
Assume $2^{\aleph_1}=\aleph_2$. \newline 
There exists a collection $\mathcal{A}=\{C_\alpha:\alpha<\omega_2\}$ of club subsets of $\omega_1$ for which $\binom{\aleph_2}{\aleph_1}\nrightarrow_{\mathcal{A}} \binom{\aleph_2}{\aleph_1}^{1,1}_2$.
\end{claim}

\par\noindent\emph{Proof}. \newline 
We commence with a general assertion which does not depend on the assumption $2^{\aleph_1}=\aleph_2$. We claim that if $\{A_\beta:\beta\in\omega_1\}$ is any collection of unbounded subsets of $\aleph_1$ then there exists a club $C\subseteq\omega_1$ such that:
\begin{enumerate}
\item [$(a)$] $\forall\beta<\omega_1, A_\beta\nsubseteq C$.
\item [$(b)$] $\forall\beta<\omega_1, A_\beta\nsubseteq \aleph_1\setminus C$.
\end{enumerate}
We construct $C$ by induction on $\varepsilon<\omega_1$. At the stage $\varepsilon=0$ we choose $a_0,b_0\in A_0$ so that $a_0<b_0$. At the stage $\varepsilon+1$ we choose $a_{\varepsilon+1},b_{\varepsilon+1}\in A_{\varepsilon+1}$ such that $b_\varepsilon<a_{\varepsilon+1}< b_{\varepsilon+1}$. If $\varepsilon$ is a limit ordinal then we let $\gamma_\varepsilon=\bigcup\limits_{\delta<\varepsilon}b_\delta$ and we choose $a_\varepsilon,b_\varepsilon\in A_\varepsilon$ such that $\gamma_\varepsilon< a_\varepsilon<b_\varepsilon$. Finally, define $C$ as the closure of $\{b_\varepsilon:\varepsilon<\omega_1\}$ in the order topology.

We first show that $\forall\beta<\omega_1, A_\beta\nsubseteq C$. Indeed, given any $\beta\in\omega_1$ we claim that $a_\beta\notin C$. For $\beta=0$, the first element of $C$ is $b_0$ and $a_0<b_0$, so $a_0\notin C$ and hence $A_0\nsubseteq C$. If $\beta=\eta+1$ then $b_\eta<a_\beta<b_\beta$ and by the construction of $C$ we can see that $C\cap(b_\eta,b_\beta)=\emptyset$ so $a_\beta\notin C$. Since $a_\beta\in A_\beta$ we infer that $A_\beta\nsubseteq C$. Similarly, if $\beta$ is a limit ordinal then $\gamma_\beta=\bigcup\limits_{\delta<\beta}b_\delta <a_\beta<b_\beta$, and $C\cap(\gamma_\beta,b_\beta)=\emptyset$ by the construction. It follows, again, that $a_\beta\notin C$ and hence $A_\beta\nsubseteq C$.
Next we show that $\forall\beta<\omega_1, A_\beta\nsubseteq \aleph_1\setminus C$. Indeed, for every $\beta\in\omega_1$ we have an element $b_\beta\in A_\beta$ which belongs to $C$ by its definition, so $b_\beta\notin\aleph_1\setminus C$ and hence $A_\beta\nsubseteq \aleph_1\setminus C$.

Let $\{A_\beta:\beta\in\omega_2\}$ enumerate all the members of $[\aleph_1]^{\aleph_1}$. Here we use the assumption $2^{\aleph_1}=\aleph_2$. By induction on $\alpha\in\omega_2$ we choose a club $C_\alpha\subseteq\omega_1$ such that $\forall\beta<\alpha, A_\beta\nsubseteq C_\alpha \wedge A_\beta\nsubseteq \aleph_1\setminus C_\alpha$. This can be done since $\{A_\beta:\beta<\alpha\}$ is a collection of $\aleph_1$ many sets.

Let $\mathcal{A}$ be $\{C_\alpha:\alpha\in\omega_2\}$. We define a coloring $c:\aleph_2\times\aleph_1\rightarrow 2$ by $c(\alpha,\beta)=0 \Leftrightarrow \beta\in C_\alpha$. Clearly, $c$ is $\mathcal{A}$-amenable. We claim that the negative relation $\binom{\aleph_2}{\aleph_1}\nrightarrow_{\mathcal{A}} \binom{\aleph_2}{\aleph_1}^{1,1}_2$ is exemplified by $c$. Indeed, if $I\in[\aleph_2]^{\aleph_2}$ and $J\in[\aleph_1]^{\aleph_1}$ then $J=A_\beta$ for some $\beta<\omega_2$. Pick up any ordinal $\alpha\in I$ so that $\beta<\alpha$. Inasmuch as $J=A_\beta\nsubseteq C_\alpha \wedge J\nsubseteq \aleph_1\setminus C_\alpha$ we conclude that $c\upharpoonright(I\times J)$ is not monochromatic. But $I,J$ were arbitrary, so we are done.

\hfill \qedref{ggggch}

\begin{remark}
\label{rrr} We make the following comments:
\begin{enumerate}
\item [$(\alpha)$] The above claim works equally well for every infinite cardinal $\kappa$ with respect to $\kappa^+$ and $\kappa^{++}$. The pertinent assumption would be $2^{\kappa^+}=\kappa^{++}$.
\item [$(\beta)$] The choice of club sets is just one example, and the method seems flexible enough to allow more instances of amenability.
\item [$(\gamma)$] The construction is taken from \cite{temp}, with little modifications. A stronger theorem is proved there under the PFA, namely there exists a collection of $\omega_2$-many club subsets of $\omega_1$ such that the intersection of any sub-collection of size $\aleph_2$ of them is finite. This might give stronger negative relations in our context.
\item [$(\delta)$] If $\mathcal{A}\subseteq[\kappa]^\kappa$ and $|\mathcal{A}|\leq\kappa$ then $\binom{\kappa^+}{\kappa}\rightarrow_{\mathcal{A}} \binom{\kappa^+}{\kappa}^{1,1}_2$ is virtually true, so we always concentrate on large enough families of $[\kappa]^\kappa$ with respect to amenable colorings.
\end{enumerate}
\end{remark}

\hfill \qedref{rrr}

The opposite direction is the content of the following:

\begin{theorem}
\label{mt} Positive relation for $\aleph_2$. \newline 
It is consistent that, for every $\mathcal{A}\subseteq[\omega_1]^{\omega_1}$ which is strongly closed under countable intersections, $\binom{\aleph_2}{\aleph_1}\rightarrow_{\mathcal{A}} \binom{\aleph_2}{\aleph_1}^{1,1}_2$ holds.
\end{theorem}

\par\noindent\emph{Proof}. \newline 
We begin by forcing the generalized Martin's axiom (Theorem \ref{sssshelah}), so $2^{\aleph_0}=\aleph_1$ and $2^{\aleph_1}>\aleph_2$.
Suppose $\mathcal{A}\subseteq[\omega_1]^{\omega_1}$ is strongly closed under countable intersections, and let $c:\aleph_2\times\aleph_1\rightarrow 2$ be any $\mathcal{A}$-amenable coloring. For every $\alpha<\aleph_2$ set $A_\alpha=\{\beta\in\omega_1: c(\alpha,\beta)=0\}$. By $\mathcal{A}$-amenability there is some $B_\alpha\in\mathcal{A}$ so that $(B_\alpha\subseteq A_\alpha)\vee(B_\alpha\subseteq \omega_1\setminus A_\alpha)$. As all we need is just $\aleph_2$-many sets from $\mathcal{A}$, we may assume without loss of generality that $B_\alpha\subseteq A_\alpha$ for every $\alpha<\aleph_2$. 

We define a forcing notion $\mathbb{P}$. A condition $(A,s)\in\mathbb{P}$ consists of $A\in[\omega_2]^{\aleph_0}$ and $s\in[\omega_1]^{\aleph_0}$. For the order, we say that $(A,s)\leq_{\mathbb{P}}(B,t)$ iff $A\subseteq B, s\subseteq t$ and $\alpha\in A\Rightarrow t\setminus s\subseteq B_\alpha$. Notice that the requirements of Theorem \ref{sssshelah} are met (in particular, $\mathbb{P}$ is $\aleph_1$-centered as each pair of conditions $(A,s),(B,s)$ is compatible and $\aleph_1^{\aleph_0}=\aleph_1$).

For every $\alpha<\omega_2$ let $D_\alpha=\{(A,s):\alpha\in A\}$. If $(A,s)\notin D_\alpha$ then $(A\cup\{\alpha\},s)\in D_\alpha$, and by the order definition we have $(A,s)\leq_{\mathbb{P}}(A\cup\{\alpha\},s)$ so $D_\alpha$ is dense.
For every $\beta<\omega_1$ let $E_\beta=\{(A,s):s\nsubseteq\beta\}$. If $(A,s)\notin E_\beta$ then we let $x=\bigcap\{A_\gamma:\gamma\in A\}$, and recall that $A_\gamma$ contains a member of $\mathcal{A}$. Since $\mathcal{A}$ is closed under countable intersections, moreover, the intersection is uncountable, there is an ordinal $\delta>\beta$ so that $\delta\in x$. Consequently, $(A,s)\leq_{\mathbb{P}}(A,s\cup\{\delta\})$ and we infer that $E_\beta$ is dense.

By Theorem \ref{sssshelah} there exists a generic set $G\subseteq\mathbb{P}$ so that $G\cap D_\alpha\neq\emptyset$ for every $\alpha<\omega_2$ and $G\cap E_\beta\neq\emptyset$ for each $\beta<\omega_1$. Set:

$$
H=\bigcup\{s:\exists A,(A,s)\in G\}.
$$

For every $\alpha\in\omega_2$ choose $(A_\alpha,s_\alpha)\in G$ such that $\alpha\in A_\alpha$. This can be done since $G\cap D_\alpha\neq\emptyset$. Recall that $2^{\aleph_0}=\aleph_1$, so for some $I\in[\omega_2]^{\omega_2}$ and a fixed $t\in[\omega_1]^{\aleph_0}$ we have $\alpha\in I\Rightarrow s_\alpha=t$. Set $J=H\setminus t$ and observe that the cardinality of $J$ is $\aleph_1$. By the construction, $c\upharpoonright(I\times J)=\{0\}$, so we are done.

\hfill \qedref{mt}

\newpage 

\section{Almost strong relations}

In the former section we focused on colorings which are amenable with respect to some $\mathcal{A}$.
We may ask what happens if $\mathcal{A}=[\omega_1]^{\omega_1}$, i.e. the usual polarized relation with no limitation on the colorings. It has been proved by Laver, \cite{MR673792}, under the assumption that there is a huge cardinal, that the relation $\binom{\aleph_2}{\aleph_1}\rightarrow \binom{\aleph_1}{\aleph_1}^{1,1}_{\aleph_0}$ is consistent. Laver indicates that Galvin announced that the stronger relation $\binom{\aleph_2}{\aleph_1}\rightarrow \binom{\tau}{\aleph_1}^{1,1}_{\aleph_0}$ for every $\tau<\omega_2$ can also be proved to be consistent from the same assumption. However, Galvin did not publish the proof.

Many years later, Jones \cite{MR2275863} used an unpublished result of Woodin in order to show the consistency of $\binom{\aleph_2}{\aleph_1}\rightarrow \binom{\tau}{\aleph_1}^{1,1}_{\aleph_0}$ for every $\tau<\omega_2$. The result of Woodin gives a special ideal over $\aleph_1$. It requires an instance of the rank-into-rank axiom I1, and it is strongly connected to the specific case of $\aleph_1$. Here we prove a general result in the spirit of Laver's proof, based only on the existence of a huge cardinal. Let us begin with the following:

\begin{definition}
\label{almdef} Almost strong polarized relations. \newline 
Assume $\kappa\leq\lambda$ are infinite cardinals, and $\tau<\lambda$ is an ordinal. \newline 
The relation $\binom{\lambda}{\kappa}\rightarrow \binom{\tau}{\kappa}^{1,1}_{\theta}$ means that for every coloring $c:\lambda\times\kappa\rightarrow\theta$ one can find $A\subseteq\lambda$ such that ${\rm otp}(A)=\tau$ and $B\in[\kappa]^\kappa$ for which $c\upharpoonright (A\times B)$ is constant. \newline 
The relation $\binom{\lambda}{\kappa}\rightarrow \binom{\lambda\ \tau}{\kappa\ \kappa}^{1,1}_2$ means that for every coloring $c:\lambda\times\kappa\rightarrow 2$ one can find either $A\in[\lambda]^\lambda, B\in[\kappa]^\kappa$ such that $c\upharpoonright (A\times B)=\{0\}$ or $A\subseteq\lambda, {\rm otp}(A)=\tau$ and $B\in[\kappa]^\kappa$ such that $c\upharpoonright (A\times B)=\{1\}$.
\end{definition}

The first relation is called the balanced almost strong polarized relation if it holds for every $\tau<\lambda$. 
The second relation (in the above definition) is the unbalanced version.
The consistency of the balanced relation for successors of regular cardinals is the main theorem of this section.

\begin{theorem}
\label{galthm} Almost strong relations. \newline 
Suppose $\theta=\cf(\theta)$ and there exists a huge cardinal above $\theta$. \newline 
Then one can force the relation $\binom{\theta^{++}}{\theta^+}\rightarrow \binom{\tau}{\theta^+}^{1,1}_{\theta}$ for every $\tau<\theta^{++}$, while preserving all cardinals and cofinalities in the interval $[\aleph_1,\theta]$.
\end{theorem}

\par\noindent\emph{Proof}. \newline 
By the existence of a huge cardinal one can force an ideal $\mathcal{I}$ which is $\theta^+$-complete and $(\theta^{++},\theta^{++},\theta)$-saturated over $\theta^+$, as shown in \cite{MR673792}. Thus, we may assume that there is a $\theta^+$-complete $(\theta^{++},\theta^{++},\theta)$-saturated ideal and $2^\theta=\theta^+$. Fix an ordinal $\tau<\theta^{++}$ (without loss of generality, $\theta^+<\tau$).

Suppose we are given a coloring $c:\theta^{++}\times\theta^+\rightarrow \theta$. For every $\alpha<\theta^{++}$ we choose $n(\alpha)\in\theta$ so that $x_\alpha=\{\beta\in\theta^+:c(\alpha,\beta)=n(\alpha)\}\in \mathcal{I}^+$. The existence of $x_\alpha$ follows from the completeness of the ideal. Let $x$ be $\{x_\alpha:\alpha<\theta^{++}\}$.
In order to control the order type of the big component in the monochromatic product, we choose a chain $(M_\eta:\eta\leq\tau)$ of elementary submodels of $\mathcal{H}(\chi)$ for some large enough regular cardinal $\chi$, satisfying the following properties for every $\eta\leq\tau$:
\begin{enumerate}
\item [$(\aleph)$] $|M_\eta|=\theta^+, \theta^+\cup\{\theta^+\}\subseteq M_\eta$.
\item [$(\beth)$] $\mathcal{I},c,\tau,x\in M_\eta$.
\item [$(\gimel)$] ${}^\theta M_\eta\subseteq M_\eta$.
\item [$(\daleth)$] If $\zeta<\eta\leq\tau$ then $M_\zeta\in M_\eta$.
\end{enumerate}
For every $\eta\leq\tau$ let $\sigma_\eta=\sup(M_\eta\cap\theta^{++})$. 

By the regularity of $\theta^{++}$ we may assume, without loss of generality, that $n(\alpha)=\iota$ for some fixed $\iota<\theta$ and every $\alpha<\theta^{++}$. 
This is true since we have a subset of $\theta^{++}$ of size $\theta^{++}$ for which $n(\alpha)=\iota$, and we can thin out the coloring only to this subset. A monochromatic product for the thinned-out coloring would be also monochromatic for the original coloring.
We may also assume that $\bigcap\limits_{\alpha\in \mathcal{C}}x_\alpha\in\mathcal{I}^+$ for every $\mathcal{C}\subseteq \theta^{++}$ of size $\theta$. The saturation of $\mathcal{I}$ ensures that this holds for some collection of $\theta^{++}$-many sets, and we may assume that this collection is all the $x_\alpha$-s.

Fix a bijection $h:\theta^+\rightarrow\tau$. Let $S_0$ be $S^{\theta^{++}}_{\theta^+}\setminus \sigma_\tau$, so $S_0$ is a stationary subset of $\theta^{++}$. For every $\delta\in S_0$ we shall try to define two sequences of ordinals:
\begin{enumerate}
\item [$(a)$] $\beta^\delta_0<\cdots<\beta^\delta_\gamma<\cdots<\theta^+$ for every $\gamma<\theta^+$.
\item [$(b)$] $\langle\alpha^\delta_{h(\gamma)}:\gamma<\theta^+\rangle$, a sequence of ordinals below $\delta$.
\end{enumerate}
The construction is done by induction on $\gamma$. Notice that the second sequence need not be increasing. For the first stage of $\gamma=0$ we choose $\beta^\delta_0=\min(x_\delta)$. Then we ask whether there exists an ordinal $\epsilon>\sigma_\tau,\epsilon<\delta$ for which $\beta^\delta_0\in x_\epsilon$. If the answer is yes then there exists $\epsilon\in M_{h(0)+1}\setminus M_{h(0)}$ such that $\epsilon<\delta$ and $\beta^\delta_0\in x_\epsilon$, by elementarity. So we choose any ordinal in $M_{h(0)+1}\setminus M_{h(0)}$ which satisfies these requirements, and this is $\alpha^\delta_0$. If the answer is no, then the process is terminated.

Assume now that $\gamma>0$, and let $\beta^\delta_\gamma$ be $\min(\bigcap\limits_{\gamma'<\gamma}x_{\alpha^\delta_{h(\gamma')}}\cap x_\delta\setminus \{\beta^\delta_{\gamma'}:\gamma'<\gamma\})$. This ordinal is well defined as the intersection is an element of $\mathcal{I}^+$ and we drop at most $\theta$-many ordinals from it, so the minimum is taken over a non-empty set. Next we ask whether there exists an ordinal $\epsilon<\delta, \epsilon>\sigma_\tau$ so that  $\{\beta^\delta_0,\ldots,\beta^\delta_\gamma\}\subseteq x_\epsilon$. If the answer is yes then there exists $\epsilon\in M_{h(\gamma)+1}\setminus M_{h(\gamma)}$ for which $\{\beta^\delta_0,\ldots,\beta^\delta_\gamma\}\subseteq x_\epsilon$ (here we use the fact that ${}^\theta M_\eta\subseteq M_\eta$ for each $\eta$, and the fact that $\{\beta^\delta_0,\ldots,\beta^\delta_\gamma\}$ is of size at most $\theta$), and we choose such an ordinal as $\alpha^\delta_\gamma$. Notice that $\alpha^\delta_{h(\gamma)}\neq\alpha^\delta_{h(\gamma')}$ for every $\gamma'<\gamma$. If the answer is no then the process is terminated and we try again at the next ordinal $\delta\in S_0$.

The induction process might be terminated, indeed, before accomplishing $\theta^+$ steps. However, we claim that for some $\delta\in S_0$ the induction holds along all the steps. For proving it, assume that for every $\delta\in S_0$ there exists an ordinal $\gamma=g(\delta)$ such that we cannot choose the required ordinals at stage $\gamma$. As mentioned above, the problem arises only for the choice of $\alpha^\delta_{h(\gamma)}$. Since $g$ is a regressive function on $S_0$, there is an ordinal $\gamma<\theta^+$ and a stationary set $S_1\subseteq S_0$ such that $\delta\in S_1\Rightarrow g(\delta)=\gamma$.

The cofinality of every $\delta\in S_1$ is $\theta^+$, and the cardinality of each sequence is at most $\theta$, so all sequences are bounded.
Applying Fodor's lemma once more, there exist a stationary set $S_2\subseteq S_1$ and an ordinal $\xi<\theta^{++}$ such that all the chosen sequences for $\delta\in S_2$ are bounded below $\xi$. 
Recall that $2^\theta=(\theta^+)^\theta=\theta^+$, so there are only $\theta^+$ many sequences of the form $(\beta^\delta_{\gamma'},\alpha^\delta_{h(\gamma')}: \gamma'<\gamma)$. We may choose, therefore, two elements $\delta_0,\delta_1\in S_2$ such that $\delta_0<\delta_1$ and they share the same sequence. But then $\delta_0$ gives a positive answer to the question that we ask at the stage of choosing $\alpha^{\delta_1}_{h(\gamma)}$, so the induction can go on for $\delta_1$, a contradiction.

We conclude that for some $\delta\in S_0$ we could define the above two sequences for every $\gamma<\theta^+$. Define $A=\{\alpha^\delta_{h(\gamma)}:\gamma<\theta^+\}$ and $B=\{\beta^\delta_\gamma: \gamma<\theta^+\}$. By the construction we have ${\rm otp}(A,<)=\tau$ and $c\upharpoonright(A\times B)=\{\iota\}$, so we are done.

\hfill \qedref{galthm}

\begin{remark}
\label{rreferee} The referee of the paper suggested a clever simplification to the construction of the sequences. We fix any $\delta\in S_0$, and we choose $\beta^\delta_0=\min(x_\delta)$. Now for every $\gamma<\theta^+$ we construct the sequences as follows. The inductive assumption is that $\langle\beta^\delta_{\gamma'}:\gamma'\leq\gamma\rangle$ and $\langle\alpha^\delta_{h(\gamma')}:\gamma'<\gamma\rangle$ were chosen. Let $B$ be $\{\beta^\delta_{\gamma'}:\gamma'\leq\gamma\}$. By the closure of each $M_\eta$ we have $B\in M_\eta$. Moreover, $M_\eta\models$ there are unboundedly many $\alpha<\theta^{++}$ for which $B\subseteq x_\alpha$. It follows that we can find some $\alpha\in M_{h(\gamma)+1}\setminus M_{h(\gamma)}$ and define it as $\alpha^\delta_{h(\gamma)}$. It means that we don't have to use Fodor's lemma at the end of the proof, and every $\delta\in S_0$ yields a monochromatic product.
\end{remark}

\hfill \qedref{rreferee}

The above theorem gives almost strong relations, as the order type of the first component can be any ordinal $\tau$ below $\theta^{++}$. There is, however, a conceptual discrepancy between almost strong relations and full strong relations. As mentioned in the introduction, the assumption $2^{\theta^+}=\theta^{++}$ rules out the strong relation $\binom{\theta^{++}}{\theta^+}\rightarrow \binom{\theta^{++}}{\theta^+}^{1,1}_2$. This is not the case when dealing with almost strong relations. The claim below generalizes an observation of Foreman (see Theorem 8.16 in \cite{MR2768681}):

\begin{claim}
\label{fforemanclm} The relation $\binom{\theta^{++}}{\theta^+}\rightarrow \binom{\tau}{\theta^+}^{1,1}_{\theta}$ for every $\tau<\theta^{++}$ is consistent with $2^{\theta^+}=\theta^{++}$.
\end{claim}

\par\noindent\emph{Proof}. \newline 
First we force $\binom{\theta^{++}}{\theta^+}\rightarrow \binom{\tau}{\theta^+}^{1,1}_{\theta}$ for every $\tau<\theta^{++}$. Now we proceed to the power set of $\theta^+$.
Let $\mathbb{P}$ be L\'evy$(\theta^{++},2^{\theta^+})$. Our claim is that the above relation still holds in the generic extension by the collapse.

For proving this fact, let $\name{f}$ be a name of a function from $\theta^{++}\times\theta^+$ into $\theta$. Choose a condition $p$ in $\mathbb{P}$ which forces $\name{f}$ to be a function. We shall define an increasing sequence of conditions $\langle p_j:j<\theta^{++}\rangle$, and a function $g:\theta^{++}\times\theta^+\rightarrow\theta$ so that $g$ belongs to the ground model.

We commence with $p_0=p$. Arriving at $j<\theta^{++}$ we choose $p_j$ so that $i<j\Rightarrow p_i\leq p_j$ and $\forall\alpha\leq j,\forall\beta<\theta^+, p_j\Vdash\name{f}(\alpha,\beta)=g(\alpha,\beta)$. This can be done because $p=p_0$ forces that $\name{f}$ is a function, hence if any condition $q$ extends $p$ and forces a value to $\name{f}(\alpha,\beta)$ then this value is uniqe. Now we use the completeness of our forcing (it is $\theta^{++}$-complete) in order to cover $\theta^+$-many pairs at each stage of the induction.

Since the forcing relation is definable in ${\rm V}$ we conclude that $g\in{\rm V}$, hence we can choose $A,B$ so that ${\rm otp}(A)=\tau,|B|=\theta^+$ and $g\upharpoonright(A\times B)$ is constant. Choose an ordinal $j<\theta^{++}$ such that $A\subseteq j$. By the construction, $p_j\Vdash\name{f}(\alpha,\beta)=g(\alpha,\beta)$ for all $\alpha\leq j$ and $\beta<\theta^+$, so $p_j$ forces that $\name{f}\upharpoonright(A\times B)$ is constant. However, $p\leq p_j$ and the choice of $p$ was arbitrary, so the empty condition forces that $\name{f}$ is constant on a product of the required size.

\hfill \qedref{fforemanclm}

What can be said about the strong polarized relation with respect to successor cardinals? Positive results in recent years demonstrated the importance of the splitting number for this issue. It turns out that the splitting number is relevant also for negative results. Suppose $B\in[\kappa]^\kappa$. We say that $S$ splits $B$ iff $|S\cap B|=|(\kappa-S)\cap B|=\kappa$. We say that $\mathcal{A}\subseteq[\kappa]^\kappa$ is a splitting family iff for every element $B\in[\kappa]^\kappa$ there exists some $S\in\mathcal{A}$ such that $S$ splits $B$. In the case of successor cardinals $\kappa=\theta^+$, there is always a splitting family over $\kappa$ of size $\kappa^+$. We need, however, an additional property:

\begin{definition}
\label{herdef} Hereditary splitting family. \newline 
Assume $\mathcal{A}\subseteq[\kappa]^\kappa$. \newline 
We call $\mathcal{A}$ a hereditary splitting family iff $\mathcal{B}\subseteq\mathcal{A}$ is a splitting family whenever $|\mathcal{B}|=|\mathcal{A}|$.
\end{definition}

The following connects hereditary splitting with negative strong relations:

\begin{theorem}
\label{mmt} Assume $\kappa<\mu=\cf(\mu)$. \newline 
If there exists a hereditary splitting family in $[\kappa]^\kappa$ of size $\mu$ then $\binom{\mu}{\kappa}\nrightarrow \binom{\mu}{\kappa}^{1,1}_2$. \newline 
Conversely, if $\binom{\mu}{\kappa}\nrightarrow \binom{\mu}{\kappa}^{1,1}_2$ and $\kappa=\cf(\kappa)$ then there exists a hereditary splitting family in $[\kappa]^\kappa$ of size $\mu$.
\end{theorem}

\par\noindent\emph{Proof}. \newline 
Let $\mathcal{A}=\{S_\alpha:\alpha<\mu\}$ be a hereditary splitting family, and define a coloring $c:\mu\times\kappa\rightarrow 2$ by $c(\alpha,\beta)=0$ iff $\beta\in S_\alpha$. We claim that $c$ exemplifies the negative relation $\binom{\mu}{\kappa}\nrightarrow \binom{\mu}{\kappa}^{1,1}_2$.

Assume towards contradiction that $c\upharpoonright(A\times B)$ is constant for some $A\in[\mu]^\mu,B\in[\kappa]^\kappa$. If $c\upharpoonright(A\times B)=\{0\}$ then $B\subseteq S_\alpha$ for every $\alpha\in A$. Consequently, the sub-collection $\mathcal{B}=\{S_\alpha:\alpha\in A\}$ is not a splitting family in $[\kappa]^\kappa$, contradicting the hereditariness assumption. Similarly, if $c\upharpoonright(A\times B)=\{1\}$ then $B\subseteq\kappa-S_\alpha$ for every $\alpha\in A$ and the same $\mathcal{B}$ is non-splitting, a contradiction.

For the opposite direction, let $c$ be a coloring which exemplifies the negative relation $\binom{\mu}{\kappa}\nrightarrow \binom{\mu}{\kappa}^{1,1}_2$. For every $\alpha<\mu$ let $S_\alpha$ be $\{\beta\in\kappa:c(\alpha,\beta)=0\}$. Set $\mathcal{A}=\{S_\alpha:\alpha<\mu\}$. 
We claim that $|\mathcal{A}|=\mu$. Indeed, without loss of generality $\alpha<\beta<\mu\Rightarrow S_\alpha\neq S_\beta$, since if some $S_\alpha$ appears $\mu$-many times then $\binom{\mu}{\kappa}\rightarrow \binom{\mu}{\kappa}^{1,1}_2$, so we may remove all the repetitions from $\mathcal{A}$ and still remain with a collection of size $\mu$.
We claim that $\mathcal{A}$ is a hereditary splitting family.

For proving this fact, assume $\mathcal{B}\subseteq\mathcal{A}$ and $|\mathcal{B}|=\mu$. Choose any $B\in[\kappa]^\kappa$ and let $A$ be $\{\alpha<\mu:S_\alpha\in\mathcal{B}\}$. Since $\binom{\mu}{\kappa}\nrightarrow \binom{\mu}{\kappa}^{1,1}_2$ as exemplified by $c$, we have $c\upharpoonright(A\times B)=\{0,1\}$. If $\mathcal{B}$ fails to split $B$ then $B\subseteq^* S_\alpha \vee B\subseteq^* \kappa-S_\alpha$ for every $S_\alpha\in\mathcal{B}$, so without loss of generality $B\subseteq^* S_\alpha$ for every $S_\alpha\in\mathcal{B}$.

Recall that $\kappa<\mu$ are regular cardinals, so we can assume without loss of generality that $B\subseteq S_\alpha$ for every $S_\alpha\in\mathcal{B}$. This can be done by removing a fixed initial segment of $\kappa$ from $B$ over $\mu$-many elements of $\mathcal{B}$. 
Recall that $A=\{\alpha<\mu:S_\alpha\in\mathcal{B}\}$ and notice that $c\upharpoonright(A\times B)=\{0\}$, a contradiction.

\hfill \qedref{mmt}

The above theorems invite further investigation, and we phrase several open problems. The strong relation $\binom{\mu}{\kappa}\rightarrow \binom{\mu}{\kappa}^{1,1}_\theta$ is balanced in the sense that the monochromatic product is of the same size for all colors. Likewise, the almost strong relation $\binom{\mu}{\kappa}\rightarrow \binom{\tau}{\kappa}^{1,1}_\theta$ for ordinals $\tau<\mu$ is balanced. One may wonder what happens at successor cardinals when dealing with the strongest unbalanced relation:

\begin{question}
\label{q0} Unbalanced relation for successor cardinals. \newline 
Suppose $\kappa$ is successor cardinal. Is it consistent that $\binom{\kappa^+}{\kappa}\rightarrow \binom{\kappa^+\ \tau}{\kappa\ \kappa}^{1,1}_2$ for every $\tau<\kappa^+$?
\end{question}

The second problem is motivated by the amenability result. We employed the generalization of Martin's axiom, for the case of $\aleph_2$. Higher generalizations are problematic. The following is natural:

\begin{question}
\label{q1} Amenable positive relations above $\aleph_2$. \newline 
Is it possible to prove the consistency of $\binom{\mu^+}{\mu}\rightarrow_{\mathcal{A}} \binom{\mu^+}{\mu}^{1,1}_2$ when $\mu>\aleph_1$, under the assumption that $2^\mu=\mu^+$ implies that $\binom{\mu^+}{\mu}\nrightarrow_{\mathcal{A}} \binom{\mu^+}{\mu}^{1,1}_2$?
\end{question}

Finally, the existence of the special ideal over $\theta^+$ can be proved when $\theta=\cf(\theta)$. One may wonder what happens at singular cardinals:

\begin{question}
\label{q2} Almost strong relations and singular cardinals. \newline 
Assume $\theta>\cf(\theta)$. Is it consistent that $\binom{\theta^{++}}{\theta^+}\rightarrow \binom{\tau}{\theta^+}^{1,1}_2$ for every $\tau<\theta^{++}$?
\end{question}

A possible direction will be to begin with a supercompact cardinal $\theta$ and a huge cardinal above it. The forcing of Laver is $\theta$-directed-closed, so if $\theta$ is Laver-indestructible then it remains supercompact after the forcing of Theorem \ref{lavthm}. Now we would like to add either Prikry of Magidor seuquence to $\theta$. The problem is to keep the special saturation property of the ideal over $\theta^+$, or to replace it by a weaker property which will be preserved by Prikry and Magidor forcing.

\newpage

\bibliographystyle{amsplain}
\bibliography{arlist}

\end{document}